\documentclass[11pt]{article} 

\usepackage[utf8]{inputenc} 

\usepackage{geometry} 
\usepackage{graphicx} 
\usepackage{amsmath}
\usepackage{amssymb}
\usepackage{amsthm}
\usepackage{subfigure}
\usepackage{tikz}
\usetikzlibrary{patterns}
\usepackage{authblk}

\newcommand{\vv}{\mathbf}
\newcommand{\erf}{\mathrm{erf}}
\usepackage{booktabs} 
\usepackage{array} 
\usepackage{paralist} 
\usepackage{verbatim} 
\usepackage{subfig} 

\usepackage{fancyhdr} 
\usepackage{sectsty}
\allsectionsfont{\sffamily\mdseries\upshape} 

\usepackage[nottoc,notlof,notlot]{tocbibind} 
\usepackage[titles,subfigure]{tocloft} 

\newtheorem{theorem}{Theorem}
\newtheorem{prop}{Proposition}
\newtheorem{lemma}{Lemma}
\newtheorem{definition}{Definition}
\newtheorem{corollary}{Corollary}

\title{Identifying reducible $k$-tuples of vectors with subspace-proximity sensitive hashing/filtering}
\author[1,2]{Gabriella Holden, 
}\author[3]{
Daniel Shiu,
}\author[1,2]{
Lauren Strutt}
\affil[1]{University of Bristol}
\affil[2]{Arqit Intern Program}
\affil[3]{Arqit Quantum Inc. {\tt daniel.shiu@arqit.uk}}

\begin{document}
\maketitle
\begin{abstract}
We introduce and analyse a family of hash and predicate functions that are more likely to produce collisions for small reducible configurations of vectors. These may offer practical improvements to lattice sieving for short vectors. In particular, in one asymptotic regime the family exhibits significantly different convergent behaviour than existing hash functions and predicates.
\end{abstract}


\section{Introduction}
The post-quantum algorithms Kyber, Dilithium, FALCON, Frodo, NTRU and others all rely on the difficulty of finding a short vector in a lattice and so understanding the difficulty of this mathematical problem in large dimensions is highly important for the choice of parameters for these algorithms. With present understanding, the most effective method for finding short vectors in a lattice of large dimension $d$ is believed to be ``sieving''. Sieving for short vectors in a lattice was first introduced by Ajtai, Kumar and Sivakumar (AKS) ~\cite{AKS}. In such sieve algorithms, a large number of vectors (exponential in $d$) of similar length are produced by taking small linear combinations of basis vectors. Pairs of vectors are then compared to see if adding/subtracting one from the other leads to a shorter vector than one of the vectors of the pair. The new, shorter vector then replaces the longer vector in the collection and the process is iterated. The process is assumed to continue until terminating when no shorter vectors can be found because the shortest vector is reached. Heuristics are used to justify the number of vectors required to be confident of continually finding new shorter vectors.

The initial AKS idea considered all possible pairs of vectors, but subsequent ideas improved on the run-time by only considering pairs of vectors that were ``close'' as these were the ones most likely to produce shorter vectors. Computing the distance between each pair would still take time quadratic in the number of vectors, but hashing each vector with a function that is statistically likely to preserve some information about location would allow sorting by hash value to find close pairs. Such a function is called a locality sensitive hash (LSH); see~\cite{polytope} for an example of this approach. A similar approach~\cite{filter} is a family of filtering predicate functions that produce a binary output for each input vector, with nearby vectors more likely to collectively pass a collection of predicates.

One extension of the sieving approach is to consider not just pairs of vectors, but triples, or more generally tuples of vectors for which a shorter linear combination can be found~\cite{tuple}. Finding such a linear combination is more involved than in the case of a pair of vectors where the process is essentially the first step in the Euclidean algorithm. Nevertheless there are fast, deterministic algorithms for finding the shortest linear combination of three~\cite{semaev} or four~\cite{nguyen} vectors. Indeed, because the complexity of sieving depends very weakly on the reduction step, one could check larger tuples for a shorter linear combination (e.g. by Schnorr-Euchner enumeration~\cite{schnorr}, or even a lower dimensional sieving approach). In this case, provided that a shorter linear combination were found a positive proportion of the time, the sieving process would still work. Reducible configurations of $k$-tuples of vectors are more common as we increase the size of $k$, and hence tuple sieving requires fewer initial vectors, which in turn reduces the memory complexity of sieving. In this paper we will consider the problem of finding reducible configurations of $k$-tuples, focussing in particular on the case $k=3$.

Again, to make tuple sieving more computationally efficient, a filtering process is required that can identify reducible tuples with less work than considering every possible tuple. This requires both a characterisation of a large family of reducible triples/tuples as well as a powerful but cheap process to identify such triples/tuples among a large collection of vectors based on this characterisation. Current approaches~\cite{klist, fast_tuple}  incrementally build up a $k$-tuple of vectors from a vector to a pair, triple, ..., $(k-1)$-tuple of vectors where each new vector is identified from pairwise hash collisions or pairwise filter survival.

In this paper, we investigate a more direct approach where reducible $k$-tuples are identified by $k$-way collisions on hash functions or a family of $k$ survivors of a filter. We further propose a generalisation of the angular sensitive hash used in \cite{angular, polytope} which can be used as a statistical test to identify common reducible configurations. The identification is related to the problem of finding sets of vectors all of which lie close to the same $(k-1)$-dimensional subspace. In this case, our methods may have applicability to classification/clustering in other high-dimensional data sets. Our generalisation also introduces locality sensitive hashes for the case $k=2$, that offer different asymptotic behaviour to existing hashes and may be of practical use for existing sieving implementations. Our analysis also suggests that filters can be adapted for families of $k$ survivors and a new family of filter function that may favour further investigation.

\subsection{Cross-polytope hashes and random projection hashes}
A cross-polytope hash function as given in \cite{polytope} is instantiated by selecting $d$ i.i.d. Gaussian distributed vectors $\vv r_i$. Then for each vector $\vv v$ input to the hash function we compute $\pm\vv r_i\cdot \vv v$. We take note of the signed index $\pm i$ of the $\pm\vv r_i$ that produced the largest scalar product and take this to be the hash value i.e. 
$$H(\vv v)=(-1)^b i: b\in\{0,1\},  (-1)^b\vv r_i\cdot\vv v\ge |\vv r_j\vv v |\ \forall 1\le j\le d$$ 
(equalities of the scalar products are vanishingly unlikely, but if a tie-breaker is necessary one may take the least $i$). Heuristically, we may think of the $\pm\vv r_i$ as approximations to the axes of a random cross-polytope since dimensionality means the $\vv r_i$ are approximately orthogonal and of the same length. The hash is likely to approximate the nearest vertex of the cross-polytope to the vector and intuitively we expect close vectors to have a common nearest vertex and so produce a hash collision.

We note some observations on the cross-polytope hash that give insight into our selection of a hash that acts as a statistical test for reducible $k$-tuples. 
\begin{itemize}
\item Firstly, we note that the cross-polytope hash tests for pairs $\vv v_1,\vv v_2$ which are close so that $\vv v_1-\vv v_2$ is likely to be shorter. It would be just as appropriate to identify pairs that are close to being diametrically opposite so that $\vv v_1+\vv v_2$ is likely to be shorter. Our intuition suggests that we should instead select the unsigned index with largest absolute value so that our hash is identical on $\pm\vv v$. 
\item Secondly, we propose that it should be possible to use fewer than $d$ random vectors and still have a test with distinguishing power. A set of $h$ random vectors would span/define a subspace of dimension $h$ and we expect that the projection of our pair of vectors into this space would roughly preserve proximity/diametrical oppositeness. Such a hash would of course be cheaper to evaluate, which could offset any loss of power. With fewer possible output values, the collision rate will increase for both true and false positive testing. This means that for smaller $h$, the family of hashes is much more amenable to empirical testing, even in large dimensions. 
\item Thirdly, we might choose to select more than a single index. If we collect the set of $a$ indices of largest absolute value, we are effectively identifying the $a$-dimensional subspace generated by vectors in our random basis that contains the largest projection of $\vv v$. Again, we anticipate that reducible pairs of vectors are more likely to share a subspace of largest projection. Alternatively, if we write $h=a+b$ an equivalent test is to collect the $b$ indices of lowest absolute value which identifies the $b$-dimensional subspace with the smallest projection. This would be the subspace most orthogonal to $\vv v$ and again we expect reducible vectors to be more likely to share this subspace of greatest orthogonality. 
\end{itemize}

For these reasons, we define a generalised family of ``random projection'' hashes as follows:

\begin{definition}
Let $\mathcal R=\{\vv r_1,\vv r_2,\ldots,\vv r_h\}\subset\mathbb R^d$ be a set of independently generated vectors whose components with respect to some orthonormal basis are distributed $\mathcal N(0,1)$. For $(a,b)$ such that $a+b=h$, we define the family of \emph{random projection hash functions} $H_{\mathcal R,a,b}(\vv v)$ according to the set $A$ of $a$ indices $i$ of the $\vv r_i$ which produce the greatest absolute scalar product with $\vv v$ in absolute size (equivalently we could record the set $B$ of size $b$ where $B=\{1,\ldots,h\}\backslash A$) i.e.
$$H_{\mathcal R,a,b}=\left\{A:\# A= a, |\vv r_i\cdot\vv v|\ge|\vv r_j\cdot\vv v|\ \  \forall i\in A,j\in \{1,\ldots,h\}\backslash A\right\}$$
ties for dot products with identical absolute size can be broken by selecting the least index.
\end{definition}

We are interested in the $k$-way collision rate of functions in this family on various $k$-tuples of input vectors. If the hash values for our $k$-tuple are uncorrelated, we anticipate a na\"\i ve collision rate of $\binom ha^{-(k-1)}$. We shall see that the actual collision rate is a function of the pairwise scalar products of our $k$-tuple. In section~\ref{sec:reducible}, we shall show that the reducibility of a $k$-tuple is also tied to these pairwise scalar products. In section~\ref{sec:empirical}, we investigate the 3-way collision rate empirically for various configurations of triples that are reducible or not quite reducible. Such empirical calculations are only feasible for relatively small values of $\binom ha^{k-1}$. In section~\ref{sec:n_analysis}, we provide an analytic expression for the collision rate in the cases where $a=1$ or $b=1$ that can be approximated numerically via a $(k+1)$-dimensional numerical integration. This is sufficient to prove a correlation between collision rates and pairwise dot products for small values of $k$ and $\min(a,b)=1$. In section~\ref{sec:asymptotic}, we consider the asymptotic behaviour of the collision rate for large $b$ with $a$ fixed and prove a generalisation of Theorem 1 of~\cite{angular}, where the 2-collision rate in the case $h=d$, $a=1$, $b=d-1$ is considered. In \cite{angular}, it is shown that the asymptotic, negative-log-probability of a collision of the signed version of the hash is $(\alpha-1)\log h$ where $\alpha(\tau):=4/(4-\tau^2)$ where $\tau$ is the distance between the two vectors. We note that this could also be written $\alpha(\vv v_1\cdot\vv v_2):=2/(1-|\vv v_1\cdot\vv v_2|)$, if we wish to relate to pairwise dot products. In the broader case we find that the negative-log-asymptotic collision rate as $b\to\infty$ can be expressed as
$$\log\frac {B((1+o(1))\alpha a,b)}{B(a,b)}$$
where $B$ is the Beta function and $\alpha$ is a function only of the pairwise scalar products of the $k$-tuple. We refer to $\alpha$ as the ``squared shortest dual diagonal'' of a set of vectors and give a mathematical definition in section~\ref{sec:asymptotic}. We also analyse the asymptotic behaviour of the collision rate for large $a$ with $b$ fixed and find different behaviour. Specifically we prove an asymptotic collision rate of
$$\Delta^{-b}\binom{a+b}b^{-k+1}$$
where $\Delta$ is the polar sine of our $k$-tuple. We also provide in passing survival rates for $k$-tuples under predicates based on the random projection hash lying above some large bound or below some small bound in absolute value. These estimates can be used to develop filtering schemes for sieving.

\section{Reducible tuples of vectors\label{sec:reducible}}
In this section we attempt to characterise families of tuples of equal length vectors\begin{footnote}{Here, we invoke a common intuition of sieving heuristics that the vectors generated by our initial process are all approximately the same length due to dimensionality and that subsequent generations of vectors are also approximately the same length.}\end{footnote} where some linear combination of the vectors is shorter than the originals.

Considering what is likely to be a more common class of reducible tuples, we first picture a typical random $k$-tuple of vectors lying on a $(d-1)$-sphere. For $k$ small relative to $d$, our $k$-tuple is likely to be close to orthogonal and thus we can picture the $2^k$ sum/differences of our $k$-tuple as the vertices of an approximate $k$-hypercube with sides of length equal to the diameter of our $(d-1)$-sphere and parallel to our vectors. We would like a minimally disruptive transformation of this hypercube, preserving approximate side length, that would guarantee a linear combination of our vectors that is shorter than the vectors themselves. A natural instinct is to compress the hypercube along a major diagonal to form a rhombic parallelepiped. If we compress to an extent that the major diagonal is less than the length of the sides, then this diagonal is a shorter linear combination. Specifically, if the vectors emanating from one vertex of the major diagonal are $\vv u_1,\vv u_2,\ldots,\vv u_k$ (and the negation of these vectors likewise converge on the opposite vertex) then the diagonal represents the vector 
$$\vv d=\vv u_1+\vv u_2+\vv+\vv u_k.$$
Each vector $\vv u_i$ is either a vector from our tuple or its negation. By expanding out $\vv d\cdot\vv d$ using the distributive law, we can see the generalisation of the cosine rule
$$||\vv d||^2=\ell^2\left(k+\sum_{i\neq j}\cos(\vv u_i,\vv u_j)\right)$$
where $\ell$ is the side length of our rhombic parallelepiped and $\cos(\mathbf u_i,\mathbf u_j)$ is the cosine of the angle between the vectors $\mathbf u_i$ and $\mathbf u_j$. To make $||\mathbf d||^2$ less than $\ell^2$ we must make the sum of the $2\binom k2$ cosine terms less than $-k+1$. The least deformative way to do this is to make each term slightly less than $-1/k$. This defines a configuration of $k$-vectors  that will serve as the boundary for our large class of reducible $k$-tuples which is where the vectors form a rhombic parallelepiped with $k$-fold rotational symmetry about its shortest diagonal which is equal to the side length. The angle between each pair of distinct vectors in this case is $\arccos(-1/k)$. Consistent with the analysis of~\cite{tuple, klist, fast_tuple} the density of random tuples that are reducible is concentrated near to this configuration and we aim to design a hash capable of testing for such configurations.

Considering the case $k=3$, our interest is drawn to triples of vectors that define oblate parallelepipeds, which is akin to triples of vectors that are close to being coplanar. Again this fits with our intuition: coplanarity is equivalent to a linear dependency in a vector space which we might think is closely related to a small integer linear combination in the corresponding lattice. More generally we can see that we are interested in $k$-dimensional parallelepipeds that have a short diagonal so that the $k$ vectors are close to lying in a $k-1$ dimensional subspace. This suggests that collisions in our family of random projection hashes could be employed to detect our tuples: the $a$ random vectors with large scalar products with our tuple vectors should all have a large component in the $(k-1)$-dimensional subspace and the $b$ other vectors should all have a large component in the orthogonal complement. Hence, our belief is that for reducible $k$-tuples $\vv v_1,\vv v_2,\ldots,\vv v_k$ we are more likely to observe a $k$-way collision for a random projection hash function $H$ than for a generic $k$-tuple.

\section{Empirical $3$-way collision rates\label{sec:empirical}}
If we restrict ourselves to linear combinations of our $k$ vectors where the coefficients are $\pm 1$, we noted in section~\ref{sec:reducible} that the reducibility of a configuration of $k$ unit vectors depends on their pairwise dot products: specifically, they are reducible if $\sum_{i\neq j}\vv u_i\vv u_j<-k-1$ for some set $\vv u_i$ of signed $\vv v_i$. Ideally our collision rate would depend only on this sum, but that turns out not to be the case. The rotational symmetry of the hash does lead us to believe that the collision rate should depend at most on the set of pairwise scalar products and we shall see that this is the case under all three investigation methods. 

For a given dimension $d$, we can generate a random $k$-tuple of vectors with prescribed pairwise scalar products by generating an arbitrary random triple, computing the Gram-Schmidt basis of the $k$-dimensional space spanned by these and then rescaling the components of each of our random vectors in the Gram-Schmidt basis. For such triples we can empirically estimate the 3-way collision rate of our random projection hash family. Na\"\i vely we expect a hash $H_{\mathcal R,a,b}$ to have a collision rate of $p\approx\binom{a+b}b^{-k+1}$ and if we sample $N$ hashes, we expect the number of collisions to be distributed $\mathcal B(N,p)$. To estimate the mean collision rate to within $\pm 1\%$ with 95\% confidence we might therefore take $N\approx -\log 0.025\times3\times 10^4\times p^{-2}\approx 10^5p^{-2}$ per a 2-sided Chernoff bound. This rules out extensive investigation of collisions with small $p$, but for a handful of $k$ and $h$ values, one can confirm that the probabilities do seem to depend at most on the pairwise scalar products and the estimates stabilise as one might expect. We would also like to experiment in a dimension large enough to demonstrate the general principle without overburdening our computation.

\subsection{Empirical results}
We investigate the case $k=3$ and consider triples whose pairwise scalar products are denoted $\alpha$, $\beta$ and $\gamma$ and the sum of the six possible ordered pairings is $2\alpha+2\beta+2\gamma=\sigma$. To do this, we work with a fixed $d=20$ (experiments were fairly stable with respect to varying dimension), and use 100,000 trials for roughly 2-3 significant figures of accuracy. We plot the data for the hash family $H_{\mathcal R,1,2}$ (i.e. where the index of the largest of three random scalar products is retained). We restrict plots to data where $\alpha,\beta,\gamma<0$ (in such cases vectors would be better reduced pairwise rather than a triple).

\begin{figure}[h]
\begin{center}
    \includegraphics[scale=0.22]{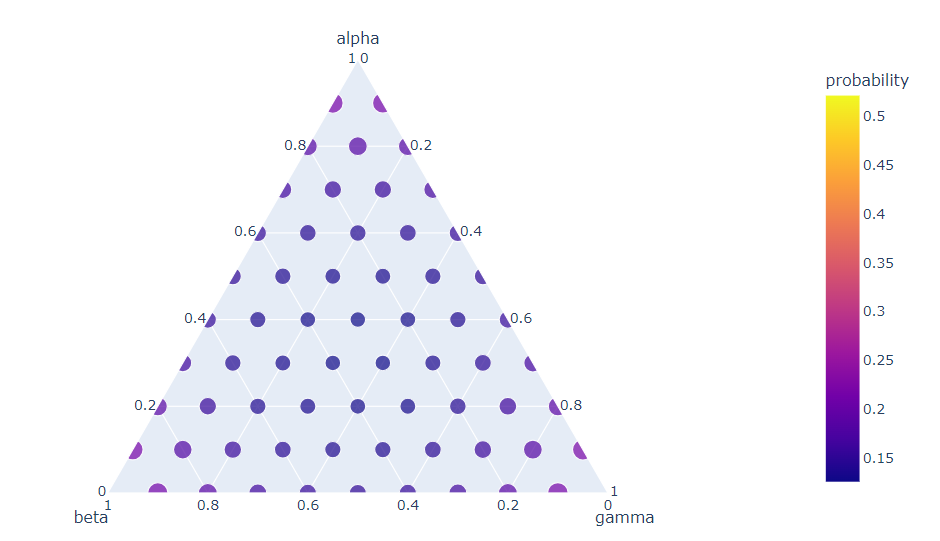}
        \includegraphics[scale=0.22]{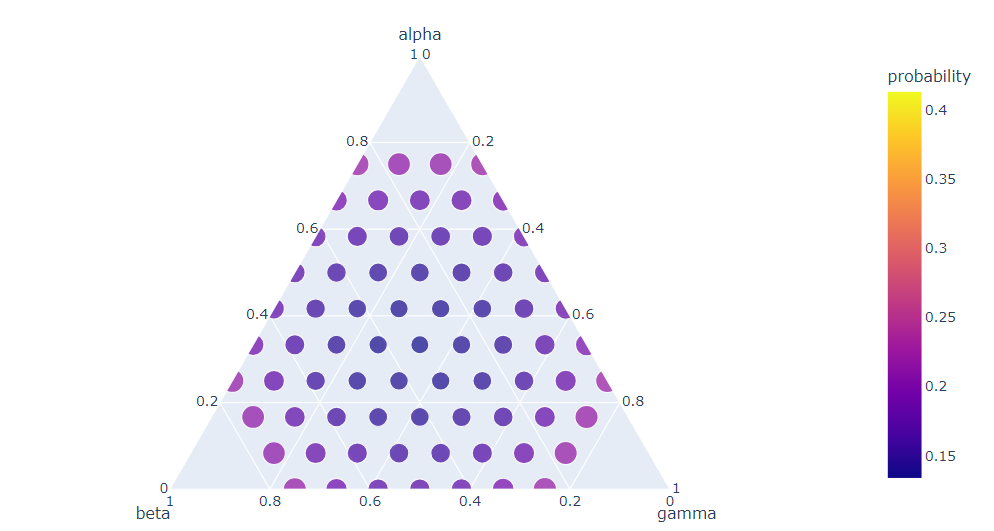}
        \includegraphics[scale=0.22]{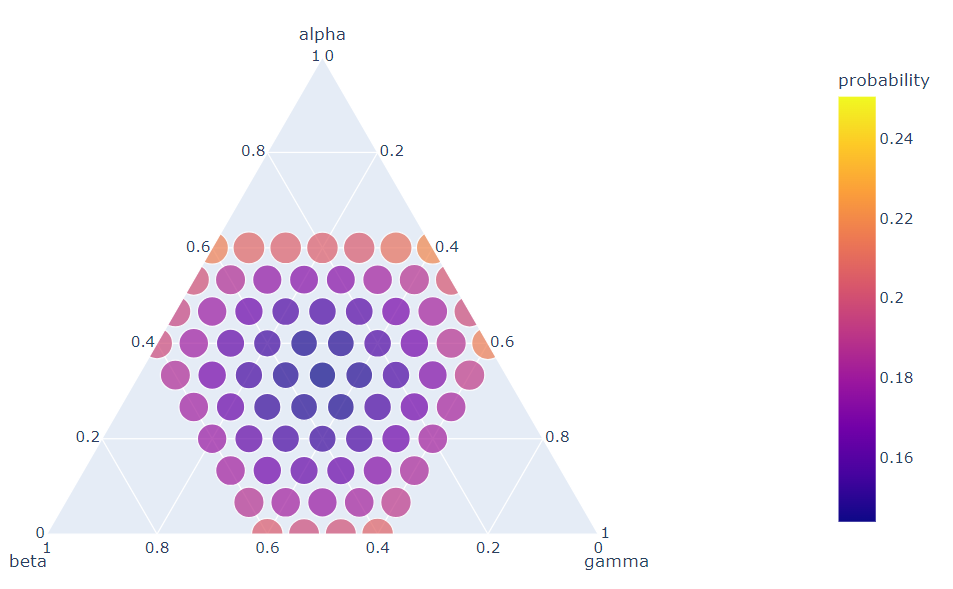}
    
    \caption{Results for $H_{\mathcal R,1,2}$ and $\sigma=-2.0, -2.4, -3.0$}
\end{center}
\end{figure}

We see the $S_3$ symmetry that we would expect in our plot and our intuition of a higher rate of collision for decreasing $\sigma$ is born out. The likelihood of a collision \emph{increases} as we move from the most common configuration which lies at the centre of our plot. For reference, the centre values in the three cases are $\approx 0.125, 0.134, 0.144$, each comfortably away from the na\~\i ve rate of $1/9$ (by roughly 12.5\%, 20.5\% and 28.9\% respectively). A good comparator for the distinguishing power of a hash collision as a test is the log-ratio of the probabilities. If we write $\rho(x,y)$ for the log-ratio of the collision rates when $\sigma=x$ and $\sigma=y$ we have
$$\rho(-2.0,-2.4)\approx 1.034,\quad \rho(-2.0,-3.0)\approx 1.035,\quad \rho(-2.0,-2.4)\approx 1.070.$$

Increasing the number of random vectors to 4, we expect the collision rates to decrease commensurate with the na\"\i ve drop from $1/9$ to $1/16$.

\begin{figure}[h]
\begin{center}
    \includegraphics[scale=0.22]{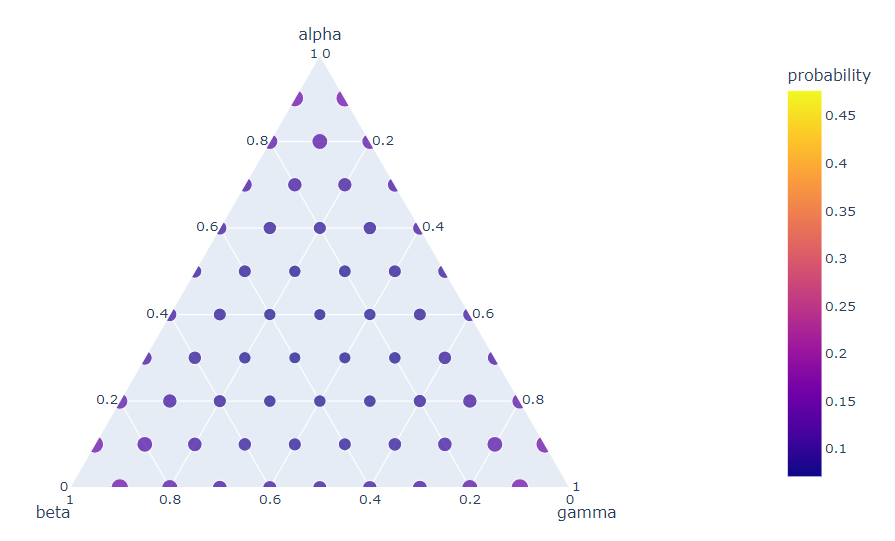}
        \includegraphics[scale=0.22]{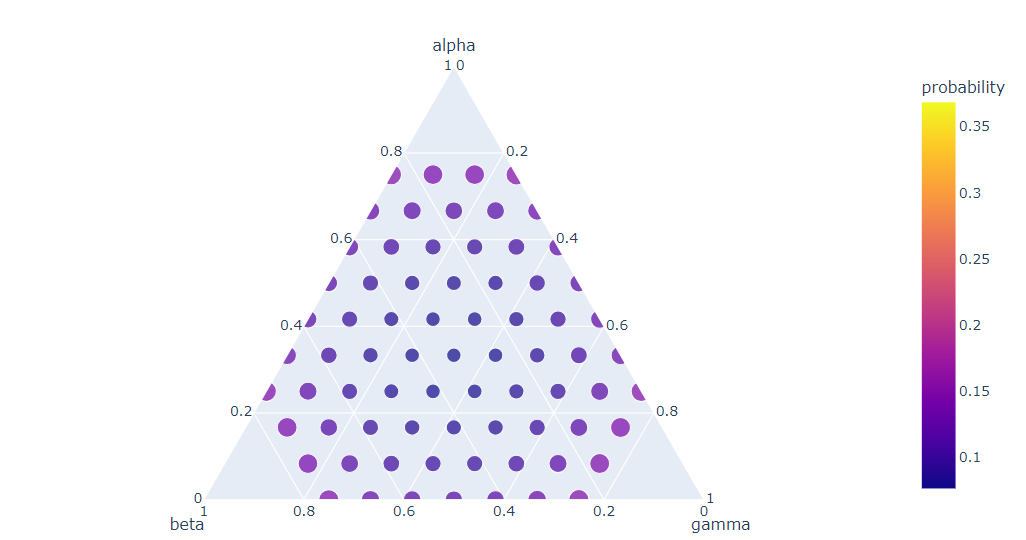}
        \includegraphics[scale=0.22]{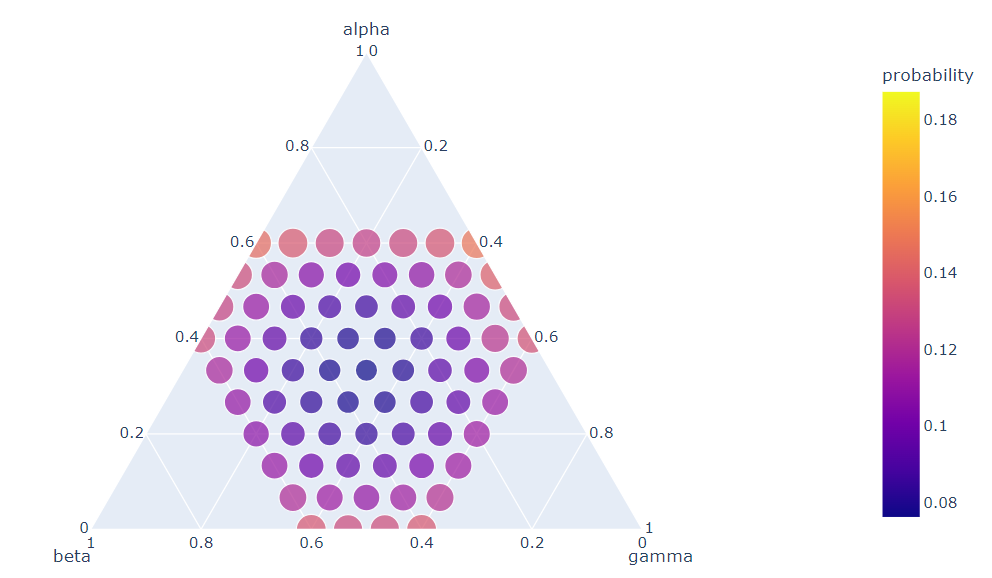}
    
    \caption{Results for $H_{\mathcal R,1,3}$ and $\sigma=-2.0, -2.4, -3.0$}
\end{center}
\end{figure}

We see a similar shape of data to the $H_{\mathcal R, 1,2}$ case. The central values are now approximately $0.0705, 0.0765, 0.0763$ (which outperform the na\"\i ve rate by roughly 12.8\%, 22.4\% and 22\% respectively). Considering log-ratios we have
$$\rho(-2.0,-2.4)\approx 1.032,\quad \rho(-2.0,-3.0)\approx 1.031,\quad \rho(-2.0,-2.4)\approx 0.999.$$
This might suggest that keeping $h$ small not only gives a cheaper hash, but also a more powerful one. However, there is clearly considerable numerical inaccuracy here.

We also consider the collision rate when we take the smallest scalar product.

\begin{figure}[h]
\begin{center}
    \includegraphics[scale=0.22]{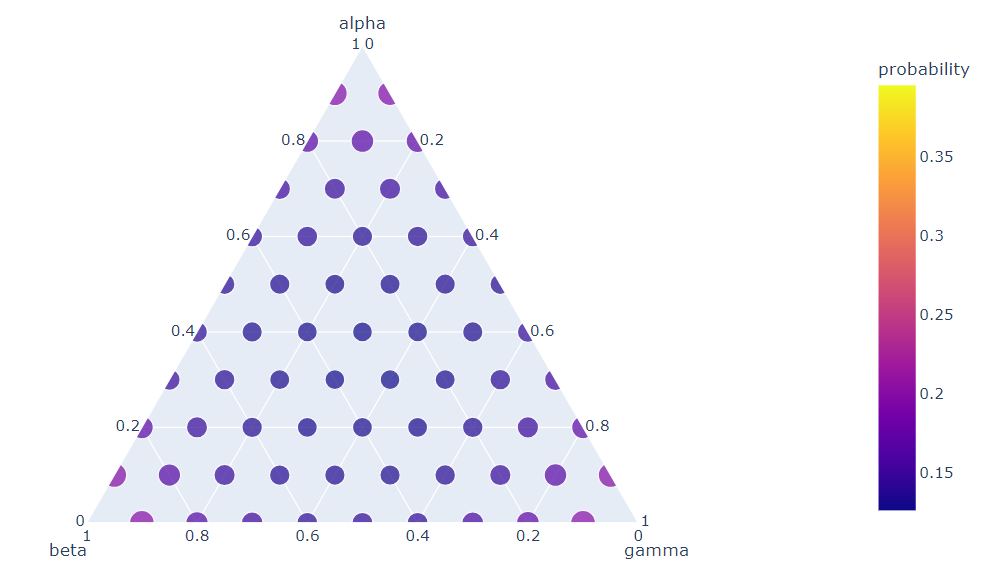}
        \includegraphics[scale=0.22]{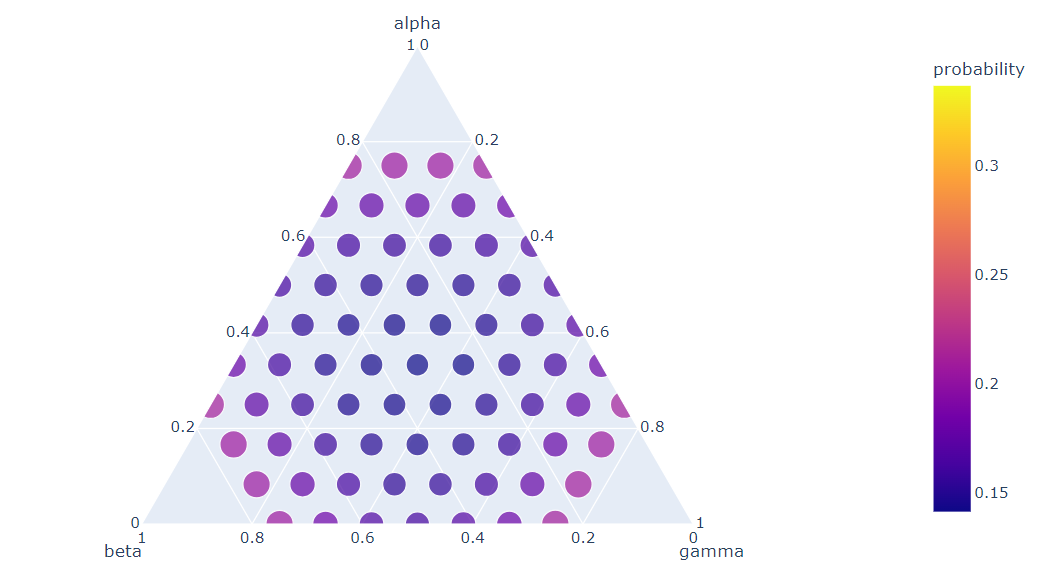}
        \includegraphics[scale=0.22]{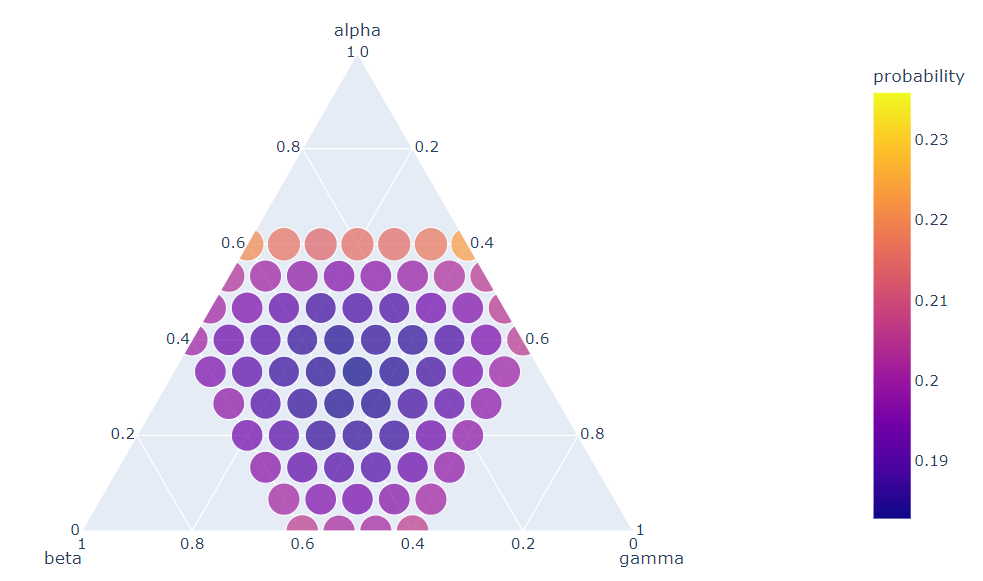}
    
    \caption{Results for $H_{\mathcal R,2,1}$ and $\sigma=-2.0, -2.4, -3.0$}
\end{center}
\end{figure}

The shape of our data does not change, but if we consider the central values we see that they are roughly $0.126, 0.141, 0.183$, exceeding the na\"\i ve rate by 13.4\%, 12.7\% and 64.5\%. Moreover, looking at log ratios we have
$$\rho(-2.0,-2.4)\approx 1.058,\quad \rho(-2.0,-3.0)\approx 1.219,\quad \rho(-2.0,-2.4)\approx 1.152.$$
This suggests that the minimal value hash is significantly stronger than the largest value hash and could benefit from deeper investigation.

Further investigating the minimum value hash with 4 random vectors we obtain

\begin{figure}[h]
\begin{center}
    \includegraphics[scale=0.22]{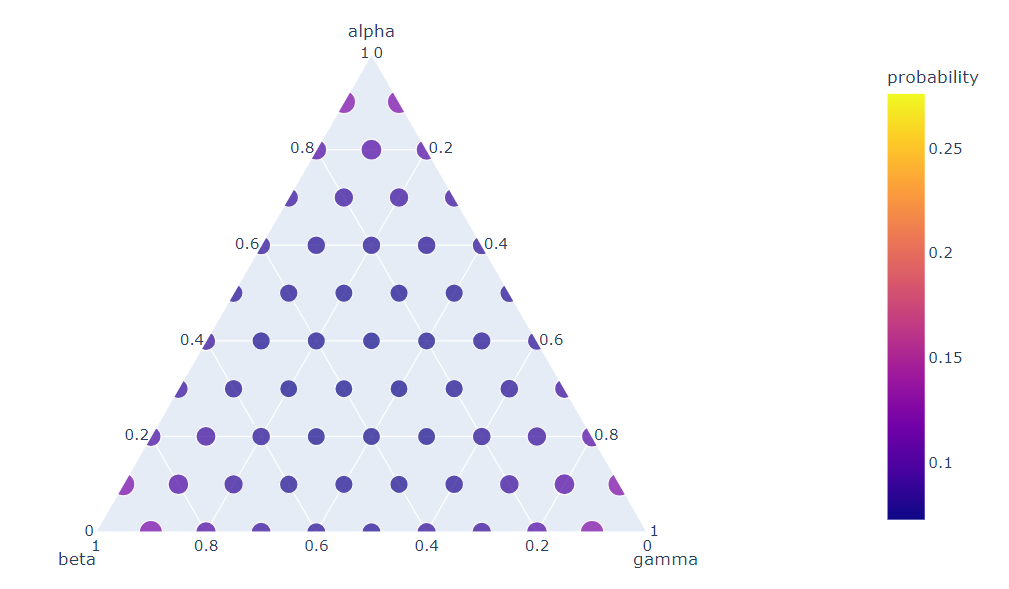}
        \includegraphics[scale=0.22]{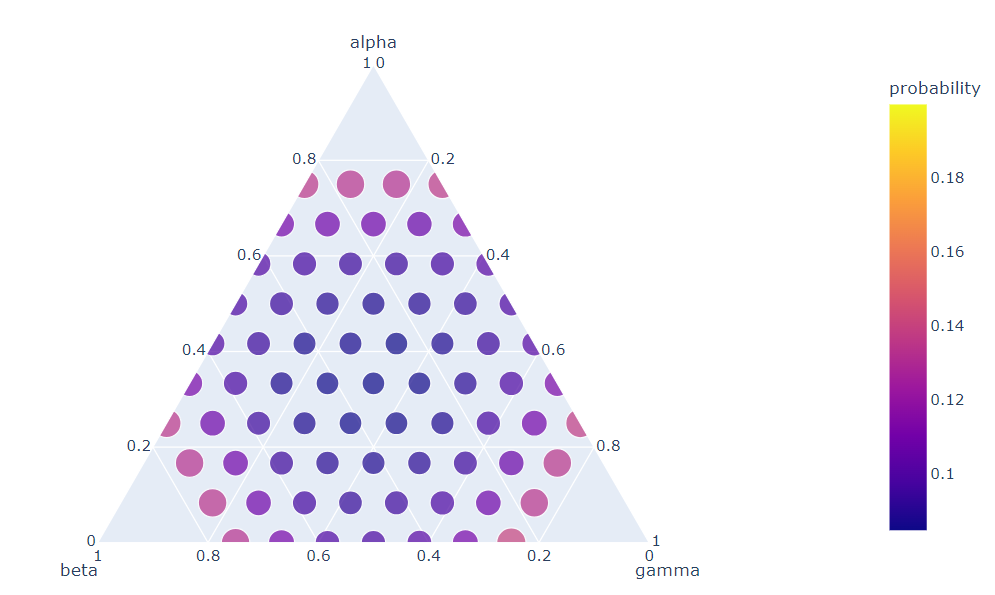}
        \includegraphics[scale=0.22]{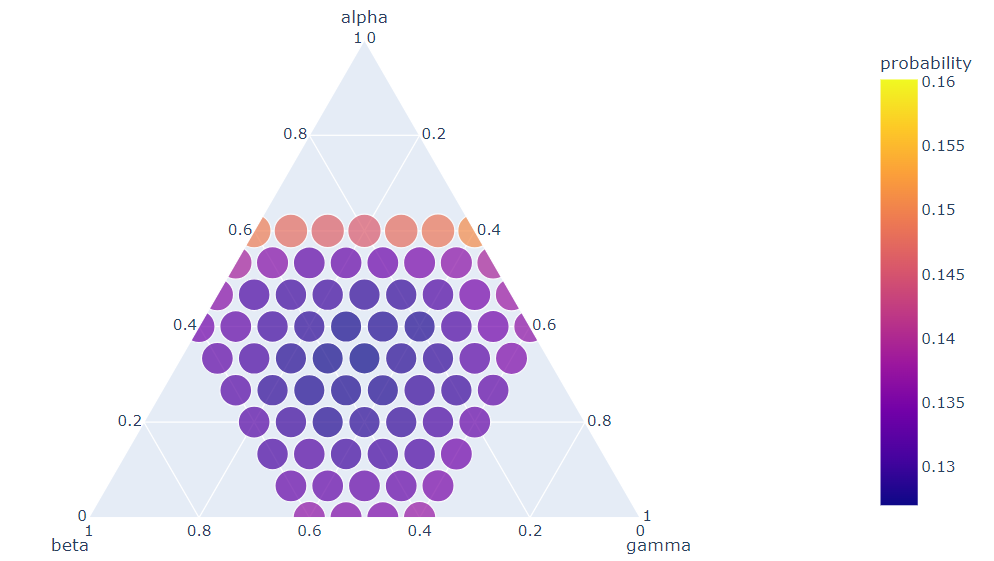}
    
    \caption{Results for $H_{\mathcal R,3,1}$ and $\sigma=-2.0, -2.4, -3.0$}
\end{center}
\end{figure}

The central values are roughly $\approx 0.0727, 0.0848, 0.127$, exceeding the na\"\i ve rate by 16.3\%, 35.7\% and 103\%. Moreover, looking at log ratios we have
$$\rho(-2.0,-2.4)\approx 1.063,\quad \rho(-2.0,-3.0)\approx 1.270,\quad \rho(-2.0,-2.4)\approx 1.196$$
so that the $H_{\mathcal R,3,1}$ hash family seems to be more powerful than the $H_{\mathcal R,2,1}$ family. Again, further investigation seems warranted.

\section{Numerically bounding the probability of a multi-collision\label{sec:n_analysis},}
In this section we develop an integral expression for the probability of a three-way collision in the case $b=1$. The generalisation of higher values of $k$ is straightforward and will be stated without proof. Given a linearly independent triple of vectors $(\vv v_1, \vv v_2,\vv v_3)$, we project into the 3-dimensional subspace spanned by these vectors. By the symmetry of the spherical distribution, the projections of the vectors $\vv r_i$ may be taken to be i.i.d. with components distributed $\mathcal N(0,1)$ with respect to an arbitrary orthonormal basis $\vv i, \vv j,\vv k$ of our subspace. By abuse of notation we also will write $\vv r_i$ for the projection of $\vv r_i$ into this subspace in the discussion below, noting that scalar products between $\vv v$s and $\vv r$s are preserved.

By symmetry, ignoring events with probability 0, 
$$\mathbb P(H(\vv v_1)=H(\vv v_2)=H(\vv v_3))=h\mathbb P(H(\vv v_1)=H(\vv v_2)=H(\vv v_3)=1).$$ 
We write $G(\vv w; \vv v_1,\vv v_2,\vv v_3)$ for the probability that for a fixed vector $\vv w$, a random Gaussian vector $\vv r$ satisfies $|\vv v_n\cdot\vv w|\le |\vv v_n\cdot\vv r|$ for $1\le n\le 3$. Then by independence,
\begin{equation}\mathbb P(H(\vv v_1)=H(\vv v_2)=H(\vv v_3)=1)=\int_{\vv w}G(\vv w; \vv v_1,\vv v_2,\vv v_3)^{h-1}d\mu(\vv w),\label{eq:pr(h=1)}\end{equation}
where $\mu$ is the Gaussian probability density measure which accounts for the probability that $\vv r_1=\vv w$.

In fact, $G$ is only a function of the pairwise scalar products of the vectors $\vv v_n$. If we write $\lambda_n=\vv w\cdot \vv v_n$, $\mu_1=\vv v_2\cdot\vv v_3$, $\mu_2=\vv v_3\cdot\vv v_1$, $\mu_3=\vv v_1\cdot\vv v_2$, then we can associate with $\vv w$ a set of 8 conjugate vectors $\{\vv w^{(c)}:\vv w^{(c)}\cdot\vv v_n=\pm \lambda_n\}$. These vectors correspond to the corners of a parallelepiped whose faces are the three pairs of parallel planes $\{\vv r:\vv r\cdot\vv v_n=\pm\lambda_n\}$ and whose edges are parallel to the vectors $\vv v_m\times\vv v_n$. The set of vectors $\vv r$ such that $|\vv v_n\cdot\vv w|\le |\vv v_n\cdot\vv r|$ for $1\le n\le 3$ is then the volume formed by the opposite solid angles at each of the 8 vertices of this parallelepiped and $G$ is the Gaussian mass of this volume.

\subsection{Calculating $G$}

\begin{center}\begin{figure}[h]\begin{tikzpicture}[scale=.45]
\draw[dotted](0,3)--(2,3);
\draw (2,3)--(9,3);
\draw[dotted](9,3)--(12,3);

\draw[dotted](0,5)--(5/3,5);
\draw[dashed] (5/3,5)--(26/3,5);
\draw[dotted] (26/3,5)--(12,5);

\draw[dotted](0,8)--(11/3,8);
\draw(11/3,8)--(32/3,8);
\draw[dotted](32/3,8)--(12,8);

\draw[dotted](0,10)--(10/3,10);
\draw (10/3,10)--(31/3,10);
\draw[dotted] (31/3,10)--(12,10);
\draw(9.4,10.5) node {$\vv w^{(c)}$};

\draw[dotted](0,0)--(5/3,5);
\draw (5/3,5)--(10/3,10);
\draw[dotted](10/3,10)--(4,12);

\draw[dotted](1,0)--(2,3);
\draw (2,3)--(11/3,8);
\draw[dotted] (11/3,8)--(5,12);

\draw[dotted](7,0)--(26/3,5);
\draw[dashed] (26/3,5)--(31/3,10);
\draw[dotted] (31/3,10)--(11,12);

\draw[dotted](8,0)--(9,3);
\draw(9,3)--(32/3,8);
\draw[dotted](32/3,8)--(12,12);

\draw[dotted](5/2,0)--(2,3);
\draw(2,3)--(5/3,5);
\draw[dotted](5/3,5)--(1/2,12);

\draw[dotted](19/2,0)--(9,3);
\draw[dashed](9,3)--(26/3,5);
\draw[dotted](26/3,5)--(15/2,12);

\draw[dotted](5,0)--(11/3,8);
\draw(11/3,8)--(10/3,10);
\draw[dotted](10/3,10)--(3,12);

\draw[dotted](12,0)--(32/3,8);
\draw(32/3,8)--(31/3,10);
\draw[dotted](31/3,10)--(10,12);

\draw[dashed, pattern=horizontal lines] (31/3,10)--(12,10)--(12,12)--(11,12);
\draw[dashed, pattern=north west lines] (31/3,10)--(10,12)--(11,12);
\draw (11.5,11.1) node {$M^{(c)}$};
\draw (6,-2) node{Diagram 1};
\end{tikzpicture}\qquad
\begin{tikzpicture}[scale=.45]
\draw (6.5,3) node {$C_1$};
\draw(2.5,9) node {$C_2$};
\draw(10,10.5) node {$C_3$};
\draw(7.3,5.8) node {$I_1$};
\draw(6,7.2) node {$I_2$};
\draw(8,7.2) node {$I_3$};
\draw(7,3) circle (4);
\draw(10,10) circle (4.5);
\draw (3,9) circle (5);
\draw[dashed] (7,3)--(10,10)--(3,9)--(7,3);
\draw[dotted] (7,3)--(7,-1);
\draw[dotted] (3,9)--(-1,12);
\draw[dotted] (10,10)--(14.5,10);
\draw(7.5,1) node{$r_1$};
\draw(1,11) node{$r_2$};
\draw(12.5,10.5) node{$r_3$};
\draw (6,-2) node{Diagram 2};
\end{tikzpicture}\end{figure}
\end{center}

The Gaussian mass of the opposite angles associated with the point $\vv w^{(c)}$ can be expressed as the integral of the Gaussian mass of the intersection of the solid angle with spheres with radius $\rho$ varying from $||\vv w^{(c)}||$ to infinity. As Gaussian density is constant on each sphere, we need only compute the proportion $\eta(\rho,\vv w^{(c)})$ of the sphere that intersects with the solid angle. Writing $M^{(c)}$ for the Gaussian mass of the opposite angle at the vertex we have
$$G=\sum_c M^{(c)};\quad\quad M^{(c)}=\frac1{\Gamma(3/2)}\int_{||\vv w^{(c)}||}^\infty \eta(\rho,\vv w^{(c)})\rho^2\exp(-\rho^2/2)d\rho.$$
The $\eta$ can be computed as the area of the intersection of three spherical caps on the sphere of radius $\rho$ divided by the area of the sphere. The three caps will be centred on $\pm\rho\vv v_1$, $\pm\rho\vv v_2$ and $\pm\rho\vv v_1$ with signs appropriate to the conjugate in question and have height $\rho-|\lambda_1|$, $\rho-|\lambda_2|$ and $\rho-|\lambda_3|$ respectively. The area of the triple cap intersection can be computed analogously to the area of the intersection of three Euclidean circles where the radii and distance between the centres of the circles are known (diagram 2). Thus the area of the intersection is the area of the triangle $I_1I_2I_3$ plus the area of the circle segments $I_1I_2$, $I_2I_3$ and $I_3I_1$. On our sphere the line segment $C_1C_2$ corresponds to a great arc of length $\rho\arccos\mu_1$ and the radius $r_1$ corresponds to the distance $\rho\arccos(\lambda_1/\rho)$ etc. For a (spherical as opposed to Euclidean) triangle with angles $A$, $B$ and $C$ opposite sides of length $a$, $b$ and $c$ respectively
$$\cos A=\frac{\cos a-\cos b\cos c}{\sin b\sin c}$$
and
$$\mathrm{Area}(\triangle abc)=A+B+C-\pi$$
using formulae from spherical trigonometry. This allows us to compute the angles $\angle I_1C_2I_3=\angle I_1C_2C_3+\angle C_1C_2I_3-\angle C_1C_2C_3$ etc. and hence the lengths $I_1I_2$ etc. and hence the area of the intersection and hence $\eta$.

\subsection{A Cartesian integral expression for the collision probability}\label{sub:cartesian}
Now that we have a means of evaluating $G$ purely in terms of the $\lambda_n$ and $\mu_n$, it makes sense to rewrite equation (\ref{eq:pr(h=1)}) as an integral over the $\lambda_n$. We use the dual basis of $\{\vv v_n\}$ as a basis for our 3-dimensional space; specifically we write $\vv w=\lambda_1\vv u_1+\lambda_2\vv u_2+\lambda_3\vv u_3$ where 
$$\vv u_1=(\vv v_2\times\vv v_3)/V, \quad \vv u_2=(\vv v_3\times\vv v_1)/V,\quad \vv u_3=(\vv v_1\times\vv v_2)/V,\quad V=\vv v_1\cdot(\vv v_2\times\vv v_3),$$
and $V$ is the scalar triple product of the $\vv v_n$. This basis is not orthogonal and the expression of the Gaussian probability density needs adjusting so that we may write it as a function of $\lambda_n$. If we write $\vv i$, $\vv j$, $\vv k$ for an arbitrary orthogonal basis of our 3-dimensional space and $\vv w=x\vv i+y\vv j+z\vv k$, then the Gaussian probability density measure is
$$d\mu(\vv w) =\frac1{(2\pi)^{3/2}}\exp(-(x^2+y^2+z^2)/2)dxdydz.$$
Now, by the properties of dual bases
$$\begin{pmatrix}\lambda_1 \\ \lambda_2 \\ \lambda_3\end{pmatrix}=B\begin{pmatrix}x\\ y\\ z\end{pmatrix}$$
where $B$ is the change of basis matrix whose rows are the vectors $\vv v_1$, $\vv v_2$, $\vv v_3$ written in the $\vv i$, $\vv j$, $\vv k$ basis (so that $\mathrm{det} B=V$). It follows that 
$$x^2+y^2+z^2=\begin{pmatrix}\lambda_1&\lambda_2&\lambda_3\end{pmatrix}(BB^T)^{-1}\begin{pmatrix}\lambda_1 \\ \lambda_2 \\ \lambda_3\end{pmatrix}$$
here $BB^T$ is the Gram matrix $M$ where $M_{ij}=\vv v_i\cdot\vv v_j$. In particular, the expression is independent of the choice of $\vv i$, $\vv j$ and $\vv k$. Writing $q_{M^{-1}}(\lambda_1,\lambda_2,\lambda_3)$ for the ternary quadratic form above, we then have
$$d\mu(\vv w) =\frac1{V(2\pi)^{3/2}}\exp(-q_{M^{-1}}(\lambda_1,\lambda_2,\lambda_3)/2)d\lambda_1d\lambda_2d\lambda_3.$$

To summarise this section, we have shown

\begin{theorem} \label{th:numerical}
Let $\mathcal V=\{\vv v_1, \vv v_2, \vv v_3\}$ be a set of three linearly independent unit vectors. Then for a random projection hash $H_{\mathcal R,h-1,1}$ for a set of $h$ vectors $\mathcal R$ chosen i.i.d. Gaussian, the probability of a three-way collision $H(\vv v_1)=H(\vv v_2)=H(\vv v_3)$ can be computed as
$$\frac h{V(2\pi)^{3/2}}\int_{-\infty}^\infty\int_{-\infty}^\infty\int_{-\infty}^\infty G_\mathcal V(\lambda_1,\lambda_2,\lambda_3)^{h-1}\exp(-q_{M^{-1}}(\lambda_1,\lambda_2,\lambda_3)/2)d\lambda_1d\lambda_2d\lambda_3$$
where $q_{M^{-1}}$ is the ternary quadratic form associated with the inverse of the Gram matrix of $\mathcal V$, $V$ is the scalar triple product of the vectors of $\mathcal V$ and $G_{\mathcal V}$ is the Gaussian mass of the exterior angles of the parallelepiped whose edges are the dual basis to $\mathcal V$ variously scaled by $\lambda_i$. In particular $G_{\mathcal V}$ can be computed as a single integral over the distance of a point from the origin.
\end{theorem}

This is the case $k=3$ of the following proposition, for which we omit the proof.

\begin{prop} 
Let $\mathcal V$ be a set of $k$ linearly independent unit vectors. Then for a random projection hash $H_{\mathcal R,h-1,1}$ for a set of $h$ vectors $\mathcal R$ chosen i.i.d. Gaussian, the probability of a $k$-way collision on the elements of $\mathcal V$ can be computed as
$$\frac h{V(2\pi)^{k/2}}\int_{-\infty}^\infty\cdots\int_{-\infty}^\infty G_\mathcal V(\lambda_1,\ldots,\lambda_k)^{h-1}\exp(-q_{M^{-1}}(\lambda_1,\ldots,\lambda_k)/2)d\lambda_1\ldots d\lambda_k$$
where $q_{M^{_1}}$ is the $k$-ary quadratic form associated with the inverse of the Gram matrix of $\mathcal V$, $V$ is the pseudoscalar of the vectors of $\mathcal V$ and $G_{\mathcal V}$ is the Gaussian mass of the exterior angles of the parallelepiped whose edges are the dual basis to $\mathcal V$ variously scaled by $\lambda_n$. In particular $G_{\mathcal V}$ can be computed as a single integral over the distance of a point from the origin.
\end{prop}

Likewise, we observe the following proposition without proof.

\begin{prop} 
Let $\mathcal V$ be a set of $k$ linearly independent unit vectors. Then for a random projection hash $H_{\mathcal R,1,h-1}$ for a set of $h$ vectors $\mathcal R$ chosen i.i.d. Gaussian, the probability of a $k$-way collision on the elements of $\mathcal V$ can be computed as
$$\frac h{V(2\pi)^{k/2}}\int_{-\infty}^\infty\cdots\int_{-\infty}^\infty F_\mathcal V(\lambda_1,\ldots,\lambda_k)^{h-1}\exp(-q_{M^{-1}}(\lambda_1,\ldots,\lambda_k)/2)d\lambda_1\ldots d\lambda_k$$
where $q_{M^{_1}}$ is the $k$-ary quadratic form associated with the inverse of the Gram matrix of $\mathcal V$, $V$ is the pseudo-scalar of the vectors of $\mathcal V$ and $F_{\mathcal V}$ is the Gaussian mass of the interior of the parallelepiped whose edges are the dual basis to $\mathcal V$ variously scaled by $\lambda_n$. In particular $F_{\mathcal V}$ can be computed as a single integral over the distance of a point from the origin.
\end{prop}

The virtue of these propositions is that they reduce an ostensibly $hk$-dimensional integral into a $k+1$ dimensional integral which is more amenable to numerical analysis techniques that provide bounds rather than probabilistic values (section \ref{sec:empirical} can be thought of as a Monte Carlo integration of such integrals). However, although highly parallelisable, these methods are still slow even for modest values of $k$. For $k=3$, it is just feasible to rigorously confirm that some the empirical estimates of section \ref{sec:empirical} are genuinely higher than the na\"\i ve collision rates as well as investigating the rate of convergence for the results of the next section for in the case of large $h$.

\section{Asymptotic collision rates}
To understand the effect of random projection hashing on the complexity of sieving algorithms, we will need an asymptotic estimate for the probabilities of $k$-way collisions and sets of $k$-survivors of predicate testing in terms of the configuration of the $k$ vectors under consideration. In this section we give a log-asymptotic estimate for the case of $H_{\mathcal R,a,b}$ with $a$ fixed and $b$ tending to infinity and an asymptotic estimate for the case of $H_{\mathcal R,a,b}$ with $b$ fixed and $a$ tending to infinity. To  state the first result, we need to clarify the dependence of the result on the configuration (pairwise scalar products) of our colliding set of vectors $V$, for which purpose we introduce the following definition.

\begin{definition}
Given a set $V$ of $k$ linearly independent unit vectors, we can refer to them as a basis for their span. Let $U=\{\vv u_1,\ldots,\vv u_k\}$ be the dual basis and consider the parallelepiped whose corners are $\pm \vv u_1\pm\vv u_2\cdots\pm\vv u_k$. We define the \emph{squared shortest dual diagonal} $\alpha(V)$ to be the square of the length of the shortest diagonal of this parallelepiped.
\end{definition}
We note that $\alpha(V)$ can be explicitly calculated from identities of vector calculus, though a simple, general closed form is elusive. We also note that in the case $k=2$ $\alpha=4/(4-\tau^2)$ per theorem 1 of \cite{angular}.

\subsection{Log-asymptotic collision rates for large $b$\label{sec:asymptotic}}
This section generalises of theorem 1 of \cite{angular} which covers the signed case $k=2$, $h=d$, $a=1$, $b=d-1$. It is worth reading the informal proof of that theorem in order to set the context for its generalisation to the case of general $k$ and sign agnosticism both when $a$ is fixed and $b$ tends to infinity and vice-versa. 

The generalisation of the proof of theorem 1 in \cite{angular} for larger values of $k$ is relatively straightforward, provided that we can prove a generalisation of Lemma 5 of \cite{angular} Appendix B. We wish to estimate the Gaussian mass of the exterior angles of a $k$-dimensional parallelepiped where the distance between the $k$-pairs of opposite faces are the same. This is a case of a tail estimate for a multivariate Gaussian distribution and although literature already exists on this topic we are content to develop our own crude log-asymptotic bound.

Throughout this section we write $\Phi_k(t)$ for the cumulative distribution function for the length of a $k$-dimensional Gaussian vector i.e.
$$\Phi_k(t)=\underset{\vv r\sim\mathcal N(0,1)^k}{\mathbb P}(|r|<t).$$
We will use the fact that for fixed $k$ as $t\to\infty$ we have $1-\Phi_k(t)=\exp(-(1+O(\log t)/t^2)t^2/2)$ and so in particular $1-\Phi_k(t)=(1-\Phi_1(t))^{1+O((\log t)/t^2)}$. 

\begin{lemma}\label{lem:exterior}
Let $G(w)$ be the Gaussian mass of the exterior angles of a $k$-dimensional parallelepiped centred on the origin, whose opposite faces are all distance $2w$ apart and whose shortest diagonal is of length $2\sqrt\alpha w$ 
$$(1-\Phi_k(w))^{\alpha +O((\log w)/w^2)}<G(w)<(1-\Phi_k(w))^{\alpha+O((\log w)/w^2)}$$
as $w\to\infty$.
\end{lemma}

\begin{proof}
We note that the region for which we wish to calculate the Gaussian mass is contained within the complement of the $k$-ball centred on the origin with radius $\sqrt\alpha w$ so that $G(w)<1-\Phi_k(\alpha w)=(1-\Phi_k(w))^{\alpha+O((\log w)/w^2)}$.

For the upper bound we consider a small parallelepiped within our region located at the corner incident on the shortest diagonal. We take side lengths $\epsilon w$ so that the furthest distance of a point in our small parallelepiped from the origin is $(1+\epsilon)\sqrt \alpha w$. Thus the Gaussian density of points within the parallelepiped is at most $\exp(-(1+\epsilon^2)\alpha w^2)$ and the volume is $(\epsilon w/2)^kV$ where $Vw^k$ is the volume of our parallelepiped. It follows that
$$G(w) > (\epsilon w/2)^kV\exp(-(1+\epsilon^2)\alpha w^2).$$
If we take $\epsilon=(\log w)/w$ we have $ (\epsilon w/2)^kV\exp(-(1+\epsilon^2)\alpha w^2)>(1-\Phi_k(w))^{\alpha+O((\log w)/w^2)}$ and the result follows.
\end{proof}

We pause to note that this lemma provides a useful asymptotic for the survival rates of $k$-tuples when filtering on the magnitude of a random projection per \cite{filter}.
\begin{corollary}
Let $P_{\vv r}(\vv v):\mathbb R^d\to\{0,1\}$ be a (random) predicate function that returns 1 if $|\vv v\cdot\vv r|>C$ and 0 otherwise for some large value $C$ where $r\sim\mathcal N(0,1)^d$. Then given a set $V\subset \mathbb R^d$ of $k$ linearly independent unit vectors
$$\mathbb P_{\vv r\sim\mathcal N(0,1)^d}(P_{\vv r}(\vv v)=1\ \forall \vv v\in V)= (1-\Phi_k(C))^{\alpha(\mathcal V)+O((\log C)/C^2)}$$
where $\alpha(\mathcal V)$ is the squared shortest dual diagonal of the vectors.
\end{corollary}

Returning to our main theme, we are now in a position to state and prove our asymptotic result

\begin{theorem}\label{th:asymptotic_B}
Suppose that $V:=\{\vv v_1, \ldots, \vv v_k\}\in S^{d-1}$ and let $\alpha(V)$ be the squared shortest dual diagonal of $V$. Let $a$ be a fixed integer and $b$ be an integer tending to infinity.  Let $H_{\mathcal R,a,b}$ be a family of random projection hashes using a set of $a+b$ random vectors whose $d$ components are i.i.d $\mathcal N(0,1)$ variables. Then
$$\mathbb P_{H\sim \mathcal H_{\mathcal R,a,b}}(H(\vv v_1)=\cdots=H(\vv v_k))=\frac {B((1+o(1))\alpha a,b)}{B(a,b)}$$
as $b\to\infty$
\end{theorem}

\begin{proof}
We fix our set of unit length vectors $V=\{\vv v_1,\ldots,\vv v_k\}$, and for any set of $a$ vectors $\mathcal R_a\subset R$ and  we define $s_{\min}(\mathcal R_a)$ to be $\min_{\vv r\in\mathcal R_a\atop1\le n\le k}|\vv v_n\cdot \vv r|$ and $\vv v(\mathcal R_a)$ to be any vector in $V$ for which this minimum is attained. Given a random Gaussian vector $\vv s$ we can see that
$$\mathbb P(|\vv s\cdot\vv v|< s_{\min}(\mathcal R_a))\le \mathbb P(|\vv s\cdot\vv v_n|<|\vv r\cdot \vv v_n|\ \forall \vv r\in\mathcal R_a,\ 1\le n\le k) \le \mathbb P(|\vv s|<s_{\min}).$$ 
By the rotational symmetry of the Gaussian distribution, the projection $\vv s$ onto $\vv v$ is distributed $\mathcal N(0,1)$ irrespective of $\vv v$ and so
$$\Phi_1(s_{\min}(\mathcal R_a))\le \mathbb P(|\vv s\cdot\vv v_n|<|\vv r\cdot \vv v_n|\ \forall \vv r\in\mathcal R_a,\ 1\le n\le k) \le \Phi_k(s_{\min}(\mathcal R_a)).$$
Moreover, for a set $\mathcal S$ of $b$ i.i.d. Gaussian vectors we see that by independence
$$\Phi_1(s_{\min}(\mathcal R_a))^b\le \mathbb P(|\vv s\cdot\vv v_n|<|\vv r\cdot \vv v_n|\ \forall \vv s\in\mathcal S,\ \vv r\in\mathcal R_a,\ 1\le n\le k) \le \Phi_k(s_{\min}(\mathcal R))^b.$$
Thus, if we write $F_s(t)$ to be the cumulative distribution function of $s_{\min}(\mathcal R_a)$ when the $R_a$ consists of $a$ i.i.d. Gaussian vectors, by the law of total probability we see that the probability $p$ of a $k$-collision in our hash function satisfies
$$\binom{a+b}a\int_0^\infty \Phi_1(t)^bdF_s(t)\le p\le\binom{a+b}a\int_0^\infty \Phi_k(t)^b dF_s(t).$$
Integrating by parts then tells us that
$$\binom{a+b}ab\int_0^\infty \Phi_1(t)^{b-1}(1-F_s(t))\Phi'_1(t)dt\le p\le\binom{a+b}ab\int_0^\infty \Phi_k(t)^{b-1}(1-F_s(t))\Phi'_k(t)dt.$$
To accommodate asymptotic estimates, we adapt the lower bounds of our integrals to tend to infinity:
$$p\ge \binom{a+b}ab\int_{t_0}^\infty \Phi_1(t)^{b-1}(1-F_s(t))\Phi'_1(t)dt,$$
$$p\le\binom{a+b}ab\int_{t_1}^\infty \Phi_k(t)^{b-1}(1-F_s(t))\Phi'_k(t)dt+\binom{a+b}a\Phi_k(t_1)^b.$$

Now, considering the tail distribution of $s_{\min}$
\begin{eqnarray*}1-F_s(t)&=&\underset{\mathcal R_a}{\mathbb P}(|\vv r\cdot\vv v_n|>t\ \forall \vv r\in\mathcal R_a,\ 1\le n\le k)\\
&=&\underset{\vv r\sim\mathcal N(0,1)^k}{\mathbb P}(|\vv r\cdot\vv v_n|>t\ ,\forall 1\le n\le k)^a \end{eqnarray*}
We write $G(V,t)$ for $\mathbb P(|\vv r\cdot\vv v_n|>t\ ,\forall 1\le n\le k)$ and note that this is the Gaussian mass of the exterior angles of the parallelepiped centred on the origin whose edges are the vectors $2t\vv v^{n}$, where $\vv v^{n}$ is the dual basis to $V$ in the vector space $\mathrm{Span}(V)$. By our lemma \ref{lem:exterior} we then have as $t\to\infty$
$$(1-\Phi_k(t))^{\alpha+O((\log t)/t^2)}<G(V,t)<(1-\Phi_k(t))^{\alpha+O((\log t)/t^2)}.$$
Thus for our upper bound we have
\begin{eqnarray*}p&\le&\binom{a+b}ab\int_{t_1}^\infty \Phi_k(t)^{b-1}(1-\Phi_k(t))^{\alpha a(1+O((\log t_1)/t_1^2))}\Phi'_k(t)dt+\binom{a+b}a\Phi_k(t_1)^b\\
&=&\binom{a+b}ab\int_{\Phi_k(t_1)}^1u^{b-1}(1-u)^{\alpha a(1+O((\log t_1)/t_1^2))}du+\binom{a+b}a\Phi_k(t_1)^b\\
&\le&\binom{a+b}a\left(bB(1+\alpha a(1+O((\log t_1)/t_1^2)),b)+\Phi_k(t_1)^b\right)\end{eqnarray*}
and for our lower bound we have
\begin{eqnarray*}p&\ge&\binom{a+b}ab\int_{t_0}^\infty \Phi_1(t)^{b-1}(1-\Phi_k(t))^{\alpha+O((\log t)/t^2)}\Phi'_1(t)dt\\
&\ge&\binom{a+b}ab\int_{t_0}^\infty \Phi_1(t)^{b-1}(1-\Phi_1(t))^{\alpha a(1+O((\log t_0)/t_0^2))}\Phi'_1(t)dt\\
&=&\binom{a+b}ab\int_{\Phi_k(t_0)}^1u^{b-1}(1-u)^{\alpha a(1+O((\log t_0)/t_0^2))}du\\
&=&\binom{a+b}abB(1+\alpha a(1+O((\log t_0)/t_0^2)),b)I_{\Phi_k(t_0)}(1+\alpha a(1+O((\log t_0)/t_0^2),b)\end{eqnarray*}
where here we use the regularised incomplete beta-function. Some beta-function identities quickly give the more presentable form
$$\binom{a+b}abB(1+\alpha a,b)=\frac{\alpha(a+b)}{\alpha a+b}\frac{B(\alpha a,b)}{B(a,b)}.$$
Choosing $t_0$ and $t_1$ to satisfy, say, $\Phi_k(t_0)=\Phi_k(t_1)=1-1/\sqrt b$ is more than sufficient to allow them to tend to infinity, whilst absorbing our auxiliary terms into our log-asymptotic as $b\to\infty$.
$$p=\frac{B(\alpha a(1+o(1)),b)}{B(a,b)}.$$\end{proof}

We note that for fixed $a$ we have $B(\alpha a(1+o(1)),b)/B(a,b)=b^{(1-\alpha)a(1+o(1))}$.

\subsection{Collision rate for large $a$}
We now prove a version of theorem \ref{th:asymptotic_B} for fixed $b$ and $a$ tending to infinity. The dominant contribution of such a probability estimate arises from the $b$th power of the Gaussian mass of the interior of a small parallelepiped centred on the origin times by the $a$th power of the Gaussian mass of its exterior angles. Close to the origin the Gaussian density is not dominated by the $\exp(-t^2/2)$ term which was useful in obliterating log-asymptotic terms in our previous proof, but is instead proportional to $t^{k-1}+O(t^{k-3})$ (which is why the introduction of $t_0$ and $t_1$ were needed in our previous proof). For our next proof, we will make use of a concentration bound around the origin. Specifically, we write $\erf$ for the Gaussian error function and note that
$$\erf(t/\sqrt 2)=\frac1{\sqrt{2\pi}}\int_{-t}^te^{-x^2/2}dx=\sqrt{\frac2\pi}t+O(t^2)$$ 
for small $t$.

\begin{lemma}\label{lem:interior_B}
Let $F(\vv s)$ be the Gaussian mass of the interior of a $k$-dimensional parallelepiped, centred on the origin, whose faces are normal to a fixed set of vectors $\vv v_1,\ldots,\vv v_k$ and are distance $2|\vv s\cdot\vv v_n|$ apart for $1\le n\le k$. Then
$$\left(1+O(||\vv s||_\Lambda^2\right)\frac1\Delta\prod_n\frac{2|\vv s\cdot\vv v_n|}{\sqrt{2\pi}}<F(\vv s)<\frac1\Delta\prod_n\frac{2|\vv s\cdot\vv v_n|}{\sqrt{2\pi}}$$
where $||\vv s||_\Lambda=\max_n|\vv s\cdot\vv v_n|$ and $\Delta$  is the volume of the parallelepiped generated by the $\vv v_n$.
\end{lemma}

\begin{proof}
Using the same coordinate system as section \ref{sub:cartesian}, our mass can be expressed as
$$F(\vv s)=\frac1{\Delta(2\pi)^{k/2}}\int_{-|\vv s\cdot\vv v_k|}^{+|\vv s\cdot\vv v_k|}\cdots\int_{-|\vv s\cdot\vv v_1|}^{+|\vv s\cdot\vv v_1|}\exp(-q_{M^{-1}}(\lambda_1,\ldots,\lambda_k)/2)d\lambda_1\ldots d\lambda_k.$$
Now if we write $q_0$ for the maximum value of $q_{M^{-1}}(\lambda_1,\ldots,\lambda_k)$ on the boundary of the hypercube $[-1,1]^k$, we have by homogeneity that 
$$0\le q_{M^{-1}}(\lambda_1,\ldots,\lambda_k)\le q_0||\vv s||_\Lambda^2$$
in our region of integration and hence
$$1- \frac{q_0||\vv s||_\Lambda^2}2\le \exp(-q_{M^{-1}}(\lambda_1,\ldots,\lambda_k)/2)\le 1$$
and
$$\left(1+O(||\vv s||_\Lambda^2)\right)\frac1\Delta\prod_n\frac{2|\vv s\cdot\vv v_n|}{\sqrt{2\pi}}|\le F(\vv s)\le\frac1\Delta\prod_n\frac{2|\vv s\cdot\vv v_n|}{\sqrt{2\pi}}$$
as required.
\end{proof}

Again, for filtering purposes, we note the following
\begin{corollary}
Let $P_{\vv r}(\vv v):\mathbb R^d\to\{0,1\}$ be a (random) predicate function that returns 1 if $|\vv v\cdot\vv r|<c$ and 0 otherwise for some small value $c$ where $r\sim\mathcal N(0,1)^d$. Then given a set $V\subset \mathbb R^d$ of $k$ linearly independent unit vectors
$$\mathbb P_{\vv r\sim\mathcal N(0,1)^d}(P_{\vv r}(\vv v)=1\ \forall \vv v\in V)= \left(1+O(c^2)\right)\frac1\Delta\left(\frac{2c^2}\pi\right)^{k/2}$$
where $\Delta$ is the volume of the parallelepiped generated by $V$.
\end{corollary}

\begin{lemma}\label{lem:exterior_B}
Let $G(\vv s)$ be the Gaussian mass of the exterior angles of a $k$-dimensional parallelepiped, centred on the origin, whose faces are normal to a fixed set of vectors $\vv v_1,\ldots,\vv v_k$ and are distance $2|\vv s\cdot\vv v_n|$ apart for $1\le n\le k$. Then
$$\left(1+O(||\vv s||_\Lambda^2\right)\prod_n\left(1-\frac{2|\vv s\cdot\vv v_n|}{\sqrt{2\pi}}\right)\le G(\vv s)\le\left(1+O(||\vv s||_\Lambda^2\right)\prod_n\left(1-\frac{2|\vv s\cdot\vv v_n|}{\sqrt{2\pi}}\right)$$
as $||\vv s||_\Lambda\to 0$ where $||\vv s||_\Lambda=\max_n|\vv s\cdot\vv v_n|$.
\end{lemma}

\begin{proof}
By the union bound, as $||\vv s||_\Lambda\to 0$
$$G(\vv s)\ge1- \sum_{n=1}^k(1-\erf(|\vv s\cdot\vv v_n|/\sqrt2))=1-\sum_{n=1}^k\sqrt{\frac2\pi}|\vv s\cdot\vv v_n|+O(||\vv s||_\Lambda^2).$$
Further, 
$$1-\sum_{n=1}^k\sqrt{\frac2\pi}|\vv s\cdot\vv v_n|=\left(1+O(||\vv s||_\Lambda^2)\right)\prod_n\left(1-\frac{2\vv |\vv s\cdot\vv v_n|}{\sqrt{2\pi}}\right).$$

Similarly, for the right hand inequality, we note that
$$\mathbb P_{\vv r\sim\mathcal N(0,1)^k}(|\vv r\cdot\vv v_m|>|\vv s\cdot\vv v_m|,|\vv r\cdot\vv v_n|>|\vv s\cdot\vv v_n|)\le(1-\erf(|\vv s\cdot\vv v_m|/\sqrt 2))(1-\erf(|\vv s\cdot\vv v_n|/\sqrt 2))$$
and then apply the second Bonferroni bound:
$$G(\vv s)\le1- \sum_{n=1}^k(1-\erf(|\vv s\cdot\vv v_n|/\sqrt 2))+\sum_m\sum_{n\neq m}(1-\erf(|\vv s\cdot\vv v_m|/\sqrt 2))(1-\erf(|\vv s\cdot\vv v_n|/\sqrt 2))$$
we conclude that
$$G(\vv s)\le\left(1+O(||\vv s||_\Lambda^2)\right)\prod_n\left(1-\frac{2\vv |\vv s\cdot\vv v_n|}{\sqrt{2\pi}}\right)$$
and the result follows.
\end{proof}

\begin{theorem}\label{th:asymptotic_A}
Suppose that $V:=\{\vv v_1, \ldots, \vv v_k\}\in S^{d-1}$. Let $b$ be a fixed integer and $a$ be an integer tending to infinity.  Let $H_{\mathcal R,a,b}$ be a family of random projection hashes using a set of $a+b$ random vectors whose $d$ components are i.i.d $\mathcal N(0,1)$ variables. Then
$$\mathbb P_{H\sim \mathcal H_{\mathcal R,a,b}}(H(\vv v_1)=\cdots=H(\vv v_k))=(1+o(1))\Delta^{-b}\binom{a+b}b^{-k+1}.$$
as $a\to\infty$, where $\Delta$ is the volume of the parallelepiped generated by the $V$.
\end{theorem}

\begin{proof}
We fix our set of unit length vectors $V=\{\vv v_1,\ldots,\vv v_k\}$, and for any finite set of vectors $\mathcal R_b$ and we define (up to events of zero probability) the vectors $\hat{\vv r}_n\in\mathcal R_b$ such that $|\vv v_n\cdot\hat{\vv r}_n|=\min_{\vv r\in\mathcal R}|\vv v_n\cdot \vv r|$ for all $1\le n\le k$. We then define
 $\vv r_{\min}(\mathcal R_b)$ to be the unique vector such that $\vv v_n\cdot\vv r_{\min}=\vv v_n\cdot \hat{\vv r}_n$ for all $1\le n\le k$. We note that $\vv r_{\min}$ does not necessarily lie in $\mathcal R_b$, though it is the corner of smallest parallelepiped whose faces are normal to the $\vv v_n$ and which contains all of $\mathcal R_b$. For a particular collision, we will require a further $a$ vectors to lie in the exterior angles of this parallelepiped. Thus similar to the proof of theorem \ref{th:asymptotic_B}, if we write $p$ for the probability of such a collision, then
$$p=\binom{a+b}a\int\cdots\int_{\mathbb R^k}G(\vv s)^adF_{\vv r_{\min}}(\vv s)$$
where $G(\vv s)$ is the Gaussian mass described in lemma \ref{lem:exterior_B}. We expect the majority of this expression's weight to be concentrated around the origin. We focus on integrating $\vv s$ over a region $\mathcal L:=\{\vv s:-\Lambda\le\vv s\cdot\vv v_n\le\Lambda, 1\le n\le k\}$ which automatically gives a lower bound and we crudely note that for $\vv s$ outside of $\mathcal L$, the mass $G$ is less than the mass of one of the $k$ regions $\{\vv w:|\vv w\cdot\vv v_n|>\Lambda\}$ which are each of mass $1-\erf(\Lambda/\sqrt2)$. It follows that
$$\idotsint_{\mathcal L}G(\vv s)^adF_{\vv r_{\min}}(\vv s)\le\frac p{\binom{a+b}a}\le\idotsint_{\mathcal L}G(\vv s)^adF_{\vv r_{\min}}(\vv s)+(1-\erf(\Lambda/\sqrt2)^a.$$
For upper and lower density functions $G_L$, $G_U$, $F_L$ and $F_U$ that satisfy $G_L\le G\le G_U$ and $F_L\le F_{\vv r_{\min}}\le F_U$ in our region $\mathcal L$, we then have
\begin{eqnarray*}\frac p{\binom{a+b}a}&\ge&\idotsint_{\mathcal L}F_{\vv r_{\min}}(\vv s)dG^a(\vv s)\\
&\ge&\idotsint_{\mathcal L}F_L(\vv s)G_L^{a-1}dG(\vv s)
\end{eqnarray*}
and likewise
$$\frac p{\binom{a+b}a}\le \idotsint_{\mathcal L}F_U(\vv s)G_U^{a-1}dG(\vv s)+(1-\erf(\Lambda/\sqrt2))^a$$
Analogous to theorem \ref{th:numerical} we change to variables $\lambda_n=\vv s\cdot\vv v_n$ for $1\le n\le k$. For our region of integration, we can take $dG=(1+O(\lambda^2))(2\pi)^{-k/2}\prod d_{\lambda_n}$ . We note that $F_{\vv r_{\min}}(\vv s)=F(\vv s)^b$ where $F(\vv s)$ is the function in lemma \ref{lem:interior_B}. Therefore writing our integral in terms of the variables $\lambda_n$, by our lemmas we may take $G_L=(1+O(\Lambda^2)) \prod_n(1-2|\lambda_n|/\sqrt{2\pi})$, $G_U=(1+O(\Lambda^2))\prod_n(1-2|\lambda_n|/\sqrt{2\pi})$, $F_L=\Delta^{-b}(1+O(\Lambda^2))^b(\prod2|\lambda_n|/\sqrt{2\pi})^b$, $F_U=\Delta^{-b}(1+O(\Lambda^2))^b(\prod2|\lambda_n|/\sqrt{2\pi})^b$. Hence
$$\frac p{\binom{a+b}a}\ge\Delta^{-b}(1+O(\Lambda^2))^{k(a+b)}\left(2a(2\pi)^{k/2}\int_0^\Lambda(1-2\lambda/\sqrt{2\pi})^{a-1}(2\lambda/\sqrt{2\pi})^bd\lambda\right)^k$$
$$\frac p{\binom{a+b}a}\ge\Delta^{-b}(1+O(\Lambda^2))^{k(a+b)}a^kB(a,b+1)^kI_{(\sqrt{2/\pi})\Lambda}(a,b+1)^k$$
where $I$ is again the normalised incomplete beta function. Likewise we have
$$\frac p{\binom{a+b}a}\le\Delta^{-b}(1+O(\Lambda^2))^{k(a+b)}a^kB(a,b+1)^k+(1-\erf(\Lambda/\sqrt2))^a.$$
Taking, say, $\Lambda=a^{-2/3}$ then gives $(1+O(\Lambda^2))^{k(a+b)}=1+O(a^{-1/3})$, $(1-\erf(\Lambda/\sqrt2)^a=O(\exp(-a^{1/3}))$ etc. The result then follows.
\end{proof}

Our theorems do not cover all instances of  the random projection hash and there may be other interesting variants to be studied by considering large $a$ and $b$ with $a/b$ tending to some (non-zero) constant.

\section{Conclusion}
We have introduced a generalisation of the cross-polytope hash/predicate family which is both practical and flexible in higher dimensional spaces. We have shown that this family of hashes (resp. predicates) produces multi-collisions (resp. multi-survivors) on reducible tuples of vectors above the na\"\i ve multi-collision rate empirically and asymptotically whilst providing tools for studying certain multi-collision rates in special cases. This should provide the means to efficiently identify reducible tuples in a non-incremental fashion. Our generalisation exhibits distinct behaviour under two different asymptotic regimes, with one regime comparable to existing hashes/predicates in the case $k=2$ (where there is log-asymptotic dependence based on the squared shortest dual diagonal of a colliding set) and the other novel (where there is asymptotic dependence on the polar sine of a colliding set).

The connection between reducibility and the squared shortest dual diagonal and between reducibility and the polar sine requires further exploration in order to develop the estimates of this paper into complexities for lattice sieving approaches, and this is a focus of continuing work.

There are further asymptotic regimes, that may exhibit other convergent behaviour for collision/survival rates and may also repay further study.

\subsection{Acknowledgements} The authors thank Daryl Burns, Roberta Faux and Sophie Stevens for many helpful discussions.

\bibliography{identifying_subspace}
\bibliographystyle{plain}
\end{document}